\documentclass[12pt]{amsart}

\newtheorem{theorem}{Theorem}[section]
\newtheorem*{theoremA}{Theorem A}
\newtheorem*{theoremB}{Theorem B}
\newtheorem{lemma}[theorem]{Lemma}

\newtheorem{corollary}[theorem]{Corollary}
\newtheorem{proposition}[theorem]{Proposition}

\theoremstyle{definition}
\newtheorem{definition}[theorem]{Definition}
\newtheorem{question}[theorem]{Question}

\theoremstyle{remark}
\newtheorem{remark}[theorem]{Remark}

\numberwithin{equation}{section}

\newcommand{\cE}{\mathcal{E}}

\newcommand{\cA}{\mathcal{A}}
\newcommand{\cF}{\mathcal{F}}

\newcommand{\cG}{\mathcal{G}}
\newcommand{\cP}{\mathcal{P}}
\newcommand{\Hh}{\mathbb{H}}
\newcommand{\D}{\mathbb{D}}
\newcommand{\C}{\mathbb{C}}
\newcommand{\Q}{\mathbb{Q}}
\newcommand{\N}{\mathbb{N}}
\newcommand{\Z}{\mathbb{Z}}
\newcommand{\St}{\mathbb{S}}
\newcommand{\R}{\mathbb{R}}

\newcommand{\de}{\delta }
\newcommand{\De}{\Delta }

\newcommand{\si}{\sigma }
\newcommand{\Si}{\Sigma }
\newcommand{\ga}{\gamma }
\newcommand{\Ga}{\Gamma }

\newcommand{\Om}{\Omega }

\newcommand{\rea}{\operatorname{Re}}

\newcommand{\ima}{\operatorname{Im}}
\newcommand{\Arg}{\operatorname{Arg}}

\newcommand{\Id}{\operatorname{Id}}

\newcommand{\bd}[1]{\partial #1}

\newcommand{\Aut}{\operatorname{Aut}}

\newcommand{\ntlim}{\operatorname{n.t.-lim}}

\begin{document}
\baselineskip=18pt

\title[Uniqueness for functional equations]{On the uniqueness of
classical semiconjugations for
self-maps of the disk}
\date{\today}
\author{Pietro Poggi-Corradini}
\address{Department of Mathematics, Cardwell Hall, Kansas State University,
Manhattan, KS 66506}
\email{pietro@math.ksu.edu}
\thanks{The author thanks Ch. Pommerenke for posing the question and
for many stimulating discussions}
\begin{abstract}
We study the uniqueness properties of classical semiconjugations for
analytic self-maps of the disk, by characterizing them as canonical
solutions to certain functional equations. As a corollary, we obtain a
complete description of all possible solutions to such equations.
\end{abstract}
\maketitle

\section{Introduction}\label{sec:intro}
Given an analytic map $\phi$ defined on the unit disk
$\D=\{z\in \C: |z|<1\}$ such that $\phi (\D)\subset\D$, and which is
not an automorphism of $\D$, we will call it a {\sf self-map} of
$\D$. Often it will be
convenient to change variables from $\D$ to the upper half-plane
$\Hh=\{z\in \C: \ima z>0 \}$ while keeping the same symbol $\phi$.
We write $\phi_{n}$ for the $n$-th iterate of $\phi$.
The pair $(\phi ,\D)$ could be thought to represent the semi-group
(under composition) $\{\phi_{n} \}_{n=0}^{\infty}$, where $\phi_{0}=\Id_{\D}$.

When $\tau$ is an automorphism of $\D$, write $\tau \in\Aut (\D)$, the
dynamics $\{\tau_{n}\}$ can be understood completely by conjugating $\tau$
to rotations, translations, or dilations. This dynamics
can be used as a model for arbitrary self-maps $\phi$ by finding
semiconjugations. In certain cases, the automorphisms of
$\D$ are not the right models and one must instead consider
automorphisms of $\C$ (the type problem for simply connected Riemann
surfaces is at the root of this fact, see \cite{cowen:1981tams}). The
present paper is mainly concerned with 
uniqueness and canonicity of the classical semiconjugations.

\subsection{The classification}\label{ssec:class}
Recall that the automorphisms of $\D$
are divided into three categories: 
\begin{enumerate}
\item {\sf elliptic automorphisms}, have one fixed point in $\D$ and
can be conjugated to rotations $z\mapsto e^{i\theta}z$ in $\D$.
\item {\sf hyperbolic automorphisms}, have no fixed points in $\D$,
two fixed points on $\bd\D$, 
and can be conjugated to dilations $z\mapsto Tz$, $T>1$, on $\Hh$
(sometimes it's more convenient to conjugate to the dilations
$z\mapsto tz$, $0<t<1$).
\item {\sf parabolic automorphisms}, have no fixed points in $\D$, one
fixed point in $\bd \D$, and can be conjugated to translations $z\mapsto
z+1$ or $z\mapsto z-1$, on $\Hh$. 
\end{enumerate}
In what follows we will include the identity map into the
elliptic category.

Likewise if $\phi$ is a self-map of $\D$, 
we again have three categories (this follows from the
classical Denjoy-Wolff Theorem and Julia's Lemma, see
\cite{bracci-pc:2003jyv} and references therein):
\begin{enumerate}
\item {\sf elliptic type}: $\phi$ has one fixed point $\zeta \in \D$ and
can be conjugated to a map that fixes the origin.
\item {\sf hyperbolic type}: $\phi$ has no fixed points in $\D$,
and has a fixed point $\zeta$ on $\bd\D$ (i.e. $\ntlim_{z\rightarrow
\zeta}\phi (z)=\zeta$) such that $\ntlim_{z\rightarrow
\zeta}\phi^{\prime} (z)=c$ exists and $0<c<1$.
Such $\phi$ can be conjugated to $z\mapsto Az+p (z)$ on $\Hh$, with
$A=1/c>1$, $\ima p (z)>0$ and $\ntlim_{z\rightarrow \infty}p (z)/z=0$
(so that $\zeta$ corresponds to $\infty$). 
\item {\sf parabolic type}: $\phi$ has no fixed points in $\D$, and
has one fixed point $\zeta$ on $\bd \D$, 
such that $\ntlim_{z\rightarrow\zeta}\phi^{\prime} (z)=1$. Such $\phi $
can be conjugated to $z\mapsto z+p (z)$ on $\Hh$, with
$\ima p (z)>0$ and $\ntlim_{z\rightarrow \infty}p (z)/z=0$.
\end{enumerate}
In all three cases the point $\zeta$ is referred to as the {\sf
Denjoy-Wolff point} of $\phi$, and $\phi_{n}$
converges locally uniformly to $\zeta$. In particular, given any point $z$ the
forward orbit $\{\phi_{n} (z) \}_{n=1}^{\infty}$ tends to $\zeta$.

The elliptic and parabolic types branch out further into subcases.
In the elliptic case, if $\phi (0)=0$ and $\phi^{\prime} (0)=0$, then
$\phi$ is {\sf 
superattracting}, while if $\phi^{\prime} (0)\neq 0$, then $\phi$ is
{\sf attracting} and in this case $|\phi^{\prime} (0)|<1$ (sometimes
this case is also called {\sf loxodromic}).
In the parabolic case, fix a point $z\in \Hh$ and consider the forward orbit
$z_{n}=\phi_{n} (z)$. Then the hyperbolic step $s_{n}:=\rho_{\Hh}
(z_{n},z_{n+1})$ decreases by Schwarz's Lemma, hence tends to a limit
$s_{\infty} (z)\geq 0$. It is  easy to see that if $s_{\infty} (z)=0$
for one point $z$, then it is zero for all $z\in \Hh$. Thus a
parabolic map $\phi$ is called {\sf zero-step} if $s_{\infty}\equiv
0$, and it is called {\sf non-zero-step} if $s_{\infty}>0$ on $\Hh$.

\subsection{The classical semiconjugations}\label{ssec:classical}
In this paper  we only consider the hyperbolic and parabolic
cases. However, it should be possible to set up the same general framework for the elliptic case and for the semiconjugations that arise in the theory of backward orbits, see \cite{pc:2000finn} and \cite{pc:2003iberoam}. 
 
One of the main methods to produce semiconjugations is to renormalize
the iterates.
Here we recall Theorem 1 of \cite{pommerenke:1979jlms}, see also
\cite{bracci-pc:2003jyv}.

Fix a self-map of $\Hh$ of hyperbolic or parabolic type, which can
therefore be written as
\begin{equation}\label{eq:phi}
\phi (z) = Az+p (z)
\end{equation}
$A\geq 1$, $\ima p (z)>0$ and $\ntlim_{z\rightarrow \infty}p (z)/z=0$
(see the classification above  in Section \ref{ssec:class}).
Write $\Aut_{\infty} (\Hh)=\{z\mapsto cz+b;\ c,b\in \R,\ c>0\}$ for
the automorphisms of $\Hh$ which fix $\infty$.
Also write $z_{n}=\phi_{n} (z_{0})=x_{n}+iy_{n}$. Then the automorphism
$\ga_{n} (z)= (z-x_{n})/y_{n}\in \Aut_{\infty} (\Hh)$ sends $z_{n}$
back to $i$. It is 
natural to consider the normalized iterates 
\begin{equation}\label{eq:gn}
g_{n}=\ga_{n}\circ \phi_n.
\end{equation}
\begin{theoremA}[Theorem 1 of \cite{pommerenke:1979jlms}]
The limit $g=\lim_{n\rightarrow \infty}g_{n}$ exists locally uniformly
in $\Hh$
and is a self-map of $\Hh$ with $g (z_{0})=i$; the sequence
$\ga_{n}\circ \ga_{n+1}^{-1}$ tends to 
$\alpha \in\Aut_{\infty} (\Hh)$; and 
\[
g\circ \phi =\alpha \circ g.
\]
Moreover, if $A>1$ in (\ref{eq:phi}) (i.e., $\phi$ is hyperbolic),
then $\alpha (z)=Az+b$ (i.e. $\alpha$ is also hyperbolic), with 
$b=\lim_{n\rightarrow \infty} (x_{n+1}-x_{n})/y_{n}\in \R$;
if $A=1$ and $\rho_{\Hh}(z_{n},z_{n+1})\downarrow s_{\infty}>0$ (i.e.,
$\phi$ is 
parabolic non-zero-step), then $\alpha (z)=z+b$ 
(i.e. $\alpha$ is parabolic), with 
$b=\lim_{n\rightarrow \infty} (x_{n+1}-x_{n})/y_{n}\neq 0$;
however, if $\phi$ is parabolic zero-step,
then $\alpha (z)=z$ and $g\equiv 1$.
\end{theoremA}
\begin{remark}\label{rem:step}
Notice that when $\phi$ is hyperbolic it can also be considered
non-zero-step. In fact, it is proved in \cite{valiron1954}, see also
Lemma 4 of \cite{bracci-pc:2003jyv}, that for every forward orbit of a
hyperbolic self-map there is a Stolz angle at infinity, i.e. a sector
$\{|\Arg z-\pi /2|<\theta \}$ with $0<\theta <\pi /2$, containing
it. Then the non-zero-step property follows from (\ref{eq:phi}).

Also in the hyperbolic case, every $b\in \R$ can arise as the constant
coefficient for $\alpha$, by changing the starting point $z_{0}$ in
the construction of $g$, see (2.9) in \cite{bracci-pc:2003jyv}.
\end{remark}
\begin{remark}\label{rem:parabnzs}
If $\phi$ is parabolic non-zero-step, and
$z_{n}=\phi_{n} (z_{0})=x_{n}+iy_{n}$, the number $b\neq 0$ in
Theorem~A is equal to $\lim_{n\rightarrow \infty}
\frac{x_{n+1}-x_{n}}{y_{n}}$. In particular, either
$\lim_{n\rightarrow \infty}x_{n}=+\infty$ or $\lim_{n\rightarrow
\infty}x_{n}=-\infty$. Also $y_{n}$ is strictly increasing, and
$x_{n}/y_{n}\rightarrow 0$, i.e., either $\Arg z_{n}\rightarrow 0$ or
$\Arg z_{n}\rightarrow \pi$ (see Remark 1 of
\cite{pommerenke:1979jlms}). In the 
following we will always assume that   $\lim_{n\rightarrow
\infty}x_{n}=+\infty$.  
\end{remark}

\subsection{Uniqueness questions}\label{ssec:uniq}
The natural questions that arise from Theorem~A are the following.
\begin{question}\label{quest:orbit}
What happens if a different orbit, other than $z_{n}=\phi_{n} (z_{0})$, is
chosen when doing the renormalization (\ref{eq:gn})?
\end{question}
\begin{question}\label{quest:uniq}
To what extent are the functions $g$ and $\alpha$ obtained in
Theorem~A unique?        
\end{question}
This last question was first raised by F. Bracci. His observation (see
also \cite{bracci-pc:2003jyv} p.~48) was that aside from Pommerenke's result,
Valiron \cite{valiron:1931bsm} had also obtain some semiconjugations using
a different 
renormalization, C. Cowen in \cite{cowen:1981tams} produced
them using the  
Uniformization Theorem (hence with a more abstract approach), and recently
Bourdon and Shapiro \cite{bourdon-shapiro:1997mams} also produced such
maps under some regularity conditions on the self-map.
Thus the question is to what extent are these maps 'essentially' the
same map.

The goal of this paper is to clarify these questions and show how the
classical semiconjugations are the canonical solutions for certain
functional equations.

\subsection{Functional equations}\label{ssec:funceq}
We are interested in the following two functional
equations:\\
{\bf Forward Equation:} {\em Given the analytic self-map $\phi$ of $\D$, find
an analytic self-map $\sigma$ of $\D$ and $\tau \in\Aut (\D)$ such
that
\begin{equation}\label{eq:forwardeq}
\sigma \circ \phi =  \tau \circ \sigma 
\end{equation}}
{\bf Backward Equation:} {\em Given the analytic self-map $\phi$ of $\D$, find
an analytic self-map $\psi$ of $\D$ and $\eta \in\Aut (\D)$ such
that
\begin{equation}\label{eq:backwardeq}
\psi \circ \eta = \phi \circ \psi    
\end{equation}}
The maps $\sigma$ and $\psi$ are called {\sf semiconjugations} and
the automorphisms $\tau$ and $\eta$ are called {\sf model automorphisms}.

It turns out, however, that an important case, the so-called ``parabolic
zero-step'' case naturally produces the following modified problem:\\
{\bf Planar Equation:} {\em Given the analytic self-map $\phi$ of $\D$, find
an analytic map $\sigma$ with $\sigma (\D)\subset\C$ and $\tau \in\Aut
(\C)$ such 
that
\begin{equation}\label{eq:planareq}
\sigma \circ \phi =  \tau \circ \sigma 
\end{equation}}
\begin{definition}\label{def:fphi}
We write $\cF (\phi)$ for the family of all non-trivial solution pairs $(\sigma
,\tau)$ to (\ref{eq:forwardeq}), i.e., we ask that $\sigma$ be
non-constant. 
\end{definition}
For instance, the pair $(g,\alpha)$ in Theorem~A is in
$\cF (\phi)$ when $\phi$ is hyperbolic or parabolic non-zero-step.
Often, it will be enough to solve these equations in a conjugation class.
In fact, if $\alpha ,\beta \in \Aut (\D)$, we have
\begin{equation}\label{eq:conjphif}
(\sigma ,\tau)\in\cF ( \phi) \Leftrightarrow (\sigma \circ \alpha, \tau)\in
\cF (\alpha^{-1} \circ \phi\circ\alpha)
\end{equation}
and
\begin{equation}\label{eq:conjtf}
(\sigma ,\tau)\in\cF ( \phi) \Leftrightarrow (\beta^{-1} \circ
\sigma,\beta^{-1} 
\circ \tau \circ \beta )\in \cF (\phi)
\end{equation}

In view of (\ref{eq:conjtf}), instead of $\cF
(\phi)$, it is enough to consider the families
\[
\cF_{e}(\phi,\theta )\cup\cF_{h}(\phi,T)\cup\cF_{p}(\phi,\pm)
\]
where
\begin{itemize}
\item $\cF_{e}(\phi,\theta )=\{\sigma:\D\rightarrow \D ;\theta\in [0,2\pi)\mid
\sigma \circ \phi =e^{i\theta}\sigma\}$.
\item   $\cF_{h}(\phi,T)=\{\sigma:\Hh \rightarrow \Hh; T>1\mid
\sigma \circ \phi 
=T\sigma\}$.
\item $\cF_{p}(\phi,\pm)=\{\sigma:\Hh \rightarrow \Hh \mid \sigma \circ \phi
=\sigma\pm 1\}$.
\end{itemize}

\begin{remark}\label{rem:conj}
With the notations of Theorem~A, when $\phi$ is
hyperbolic and $\alpha (z)=Az+b$ ($A>1$), we can let $\beta (z)=z-b/ (A-1)$ in
(\ref{eq:conjtf}), then $\beta^{-1}\circ \alpha \circ\beta
(z) =Az$, hence 
$\beta^{-1} \circ g\in \cF (\phi,A)$. On the other hand, when $\phi$ is
parabolic non-zero-step and $\alpha (z)=z+b$, we can let $\beta (z)=|b|z$ in
(\ref{eq:conjtf}), then $\beta^{-1}\circ \alpha \circ\beta
(z) =z+b/|b|$, hence $\beta^{-1} \circ g\in \cF (\phi ,\pm 1)$.
\end{remark}

Before solving these functional equations we must describe some
further properties of the classical semiconjugations.

\section{Properties of the classical semiconjugations}\label{sec:properties}

\subsection{Hyperbolic maximality}\label{ssec:metricmax}
The following observation is due to Ch.~Pommerenke (oral communication).
The semiconjugations obtained in Theorem~A are ``maximal''
with respect to the hyperbolic metric.
\begin{proposition}\label{pro:proposicion}
Suppose $\phi$ is either hyperbolic or parabolic, as in
(\ref{eq:phi}), and let $g$ be the
semiconjugation obtained in Theorem~A. If $\phi$ is
hyperbolic or parabolic non-zero-step, then for every $\sigma$ 
such that
$(\sigma,\tau)
\in\cF (\phi)$ for some $\tau$, we have
\begin{equation}\label{eq:maximal}
\rho_{\Hh} (\sigma (z),\sigma (w))\leq \rho_{\Hh} (g (z),g (w))\qquad
\forall z,w\in \Hh. 
\end{equation}
If $\phi$ is parabolic zero-step, then $\cF (\phi)$ is empty.
\end{proposition}

\begin{proof}[Proof of Proposition \ref{pro:proposicion}]
Write $\tau_{n}$ for the $n$-th iterate of $\tau$. Then, given $z,w\in \Hh$,
\begin{eqnarray*}
\rho_{\Hh} (\sigma (z),\sigma (w))
& = &
\rho_{\Hh} (\tau_{n}\circ \sigma (z),\tau_{n}\circ \sigma (w))\\
& = &
\rho_{\Hh} (\sigma\circ \phi_{n} (z),\sigma\circ \phi_{n}(w))\\
& \leq &
\rho_{\Hh} (\phi_{n} (z),\phi_{n}(w))\\
& = &
\rho_{\Hh} (g_{n} (z),g_{n}(w)),
\end{eqnarray*}
where $g_{n}$ is defined in (\ref{eq:gn}).
Proposition \ref{pro:proposicion} now follows from Theorem~A.
\end{proof}

\subsection{Univalence and covering properties of the classical
semiconjugations}\label{ssec:classuniv} 
The semiconjugations obtained in Theorem~A have
the following useful univalence and covering properties.
Suppose $\phi$ is as in (\ref{eq:phi}), and suppose $\phi$ is either
hyperbolic or parabolic non-zero-step. Let $(g,\alpha )$ be the pair
of functions obtained in Theorem~A, and let $\ga_{n}$ be
as in (\ref{eq:gn}). 
\begin{lemma}[Lemma 2 and Theorem 3 of
\cite{pommerenke:1979jlms}]\label{lem:pom}
The sequence
$\alpha_{n}^{-1}\circ g \circ \ga_{n}^{-1}$ 
converges uniformly on compact
subsets of $\Hh$ to 
the identity map on $\Hh$. 
Moreover, sequences $\rho_{k}, M_{k}\uparrow+\infty$ can be chosen so that
$g$ is univalent on the set 
\begin{equation}\label{eq:u}
U=\bigcup_{k=1}^{\infty}\bigcup_{n\geq M_{k}} (x_{n}+y_{n}\De
(i,\rho_{k}))
\end{equation}
Furthermore,
\begin{enumerate}
\item [(a)] If $\phi$ is hyperbolic, $U$ and $g (U)$ have an {\sf
inner-tangent at infinity}, i.e., given  $\delta_{n}\downarrow 0$, there is
$R_{n}\uparrow \infty$ such that $\{|z|>R_{n}, \delta_{n}<\Arg z<\pi
-\delta_{n} \}\subset\Omega$. Finally, $g$ has the
following {\sf isogonality} properties:
\begin{equation}\label{eq:isogonal}
\ntlim_{z\rightarrow \infty}g (z)=\infty \mbox{ and
}\ntlim_{z\rightarrow \infty}\frac{\Arg g (z)}{z}=0  
\end{equation}
\item [(b)] If $\phi$ is parabolic non-zero-step and $x_{n}\rightarrow
+\infty$, then the region $V=g (U)$, has a {\sf
lateral-tangent at $+\infty$ with respect to $\Hh$}, i.e., given
$\de_{n}\downarrow 0$ 
there is $t_{n}\uparrow+\infty$ such that $\{\rea z>t_{n},\
\de_{n}<\ima z<1/\de_{n}\}\subset V$.
\end{enumerate}
\end{lemma}
\begin{remark}\label{rem:rem2}
When $\phi$ is parabolic non-zero-step and $z_{n}=x_{n}+iy_{n}$ is an
orbit of $\phi$, then $y_{n}$ increases, so $y_{n}\uparrow L_{\infty}\leq
+\infty$. In particular, if $L_{\infty}<\infty$, then $U$
itself has a lateral tangent at $+\infty$. But, if $L_{\infty}=\infty$,
it is not clear. By Remark 2 of
\cite{pommerenke:1979jlms}, it turns out that $L_{\infty}<\infty$ if and
only if $g$ has a finite  angular derivative at infinity, i.e.,
\[
\ntlim_{z\rightarrow \infty}\frac{g (z)}{z}=\frac{1}{L_{\infty}}.
\]
\end{remark}
\begin{proof}{Proof of Lemma \ref{lem:pom}}
In \cite{pommerenke:1979jlms} it is proved that
$\alpha_{n}^{-1}\circ g \circ \ga_{n}^{-1}\rightarrow \Id_{\Hh}$ and
that for any $\rho >0$ there is $m_{\rho}\in \N$ so that $g$ is
univalent on 
\begin{equation}\label{eq:omrho}
\Om (\rho ,m_{\rho}) =\bigcup_{n\geq m_{\rho }} (x_{n}+y_{n}\De
(i,\rho))=\bigcup_{n\geq m_{\rho }}\ga_{n}^{-1} (\De
(i,\rho))
\end{equation}
Now, suppose a sequence $\rho_{k}\uparrow +\infty$ is given. Then
$g$ is univalent on $\Om (\rho_{1},M_{1})$. Suppose
that $M_{1}<\cdots<M_{k}$ have been chosen so that $g$ is univalent on
\[
A_{k}:=\bigcup_{j=1}^{k}\Om (\rho_{j},M_{j})
\]
Let $m_{\rho_{k+1}}$ be as in
(\ref{eq:omrho}). Choose $M_{k+1}>m_{\rho_{k+1}}$ so large that, for
all $n\geq M_{k+1}$,
\begin{equation}\label{eq:disjoint}
g\left(\ga_{n}^{-1} (\De
(i,\rho_{k+1})) \right)\cap g\left(A_{k}\setminus \Om
(\rho_{k+1},m_{\rho_{k+1}})\right)=\emptyset.  
\end{equation}
This can be done because $g\circ \ga_{n}^{-1}-nb\rightarrow
\Id_{\Hh}$. We claim that $g$ 
is univalent on $A_{k+1}$. In fact, suppose $g (z)=g (w)$ for $z,w\in
A_{k+1}$. If $z\neq w$, then necessarily, $z\in A_{k+1}\setminus
A_{k}$ and $w\in 
A_{k}\setminus \Om (\rho_{k+1},m_{\rho_{k+1}})$ (or viceversa).
Hence, $z\in\ga_{n}^{-1} (\De
(i,\rho_{k+1}))$ for some $n\geq M_{k+1}$, and by (\ref{eq:disjoint}),
$g (z)\neq g (w)$, therefore $z=w$.

Furthermore,
part (a) is well-known, see for instance p.~47 of
\cite{bracci-pc:2003jyv} and Lemma 5.1 of \cite{pc:2000finn}.

For part (b), let $(g,\alpha)$ be given by
Theorem~A, and assume
$\lim_{n\rightarrow \infty}x_{n}=+\infty$, so that $\alpha (z)=z+b$ with
$b>0$. 
Consider the half-strips 
$\St_{t}=\{\rea z>t,\
\de<\ima z<1/\de\}$, and pick
a large hyperbolic disk $\De (i,\rho)=\{z\in \Hh : \rho_{\Hh}
(z,i)<\rho \}$, so that $\cup_{k\geq 0}\left(\De (i,\rho)+kb
\right)\supset\St_{0}$.
Now consider the compact set $B=\{z\in \Hh:\ \rho_{\Hh} (z,i)\leq
2\rho\}$ and pick $\epsilon <\rho$. Then, by Lemma
\ref{lem:pom} there exists 
$N\in \N$ such that
\[
\rho_{\Hh} (\alpha_{n}^{-1}\circ g \circ \ga_{n}^{-1}
(z),z)=\rho_{\Hh} (g \circ \ga_{n}^{-1} (z),\alpha_{n} (z))<\epsilon 
\] 
for $n>N$, and uniformly for $z\in \bd B$. Hence, for
$n>\max\{N,m_{\rho}\}$, and for every $w\in
\De (i,\rho)$, we have that $g \circ \ga_{n}^{-1} (\bd B)$ winds
around $\alpha_{n} (w)$ exactly once, i.e., $\alpha_{n} (w)\in g
(U)$. 
Thus $\St_{t}\subset g(U) $ for $t$ sufficiently large.  
\end{proof}

\subsection{Orbits and canonical solutions}\label{ssec:orbits}
Given a point $z\in \D$, define $\{\phi_{n} (z) \}_{n=0}^{+\infty}$ to
be the {\sf forward orbit} starting at $z$; $\cup_{n=0}^{+\infty}
\phi_{n}^{-1} (z)$ the {\sf backward orbit} starting at $z$.
A sequence $\{w_{n} \}_{n=1}^{\infty }$ is called
a {\sf backward iteration sequence} if $\phi^{-1} (w_{n})\neq
\emptyset$ for $n=1,2,3,\dots$, and $w_{n+1}\in \phi^{-1} (w_{n})$.
 
Moreover, we say that
$\cup_{k=0}^{+\infty}\cup_{n=0}^{+\infty} \phi_{n}^{-1} (\phi_{k}
(z))$ is the {\sf grand orbit} generated by $z$. If $z,w$ belong to
a grand orbit, then there are $n,k\in \N$ so that $\phi_{n}
(z)=\phi_{k} (w)$. This can be taken to define an equivalence relation
on $\D$, and then the grand orbits are the equivalent classes.
Write $\cG (\phi)$ for the set of all grand orbits for the self-map
$\phi$. We do not worry here about what extra structures can be put on
$\cG (\phi)$.
If $\sigma$ is a semiconjugation solving (\ref{eq:forwardeq}), then it
sends grand orbits for $\phi$ into grand orbits for $\tau$, i.e.,
$\sigma$ induces a map $\Sigma$ from $\cG (\phi)$ to $\cG (\tau)$.
\begin{remark}\label{rem:fpaut}
If $p\in \D$ is a fixed point of $\phi$, then $\sigma (p)$ is a fixed
point of $\tau$.
\end{remark}
\begin{definition}\label{def:canonical}
We say that a solution $\sigma$ of (\ref{eq:forwardeq}) is {\sf
canonical} if the induced map $\Sigma$ from  $\cG (\phi)$ to $\cG
(\tau)$ is a bijection.
\end{definition}

\begin{proposition}\label{prop:gcanon}
The semiconjugation $g$ obtained in Theorem A is canonical.
\end{proposition}
\begin{proof}
Write $\langle z\rangle$ for
the grand orbit generated by $z$ and $\Si_{g}$ for the map induced by
$g$ between the spaces of grand orbits.

In the hyperbolic case, recall that by Remark \ref{rem:rem2}
every forward orbit $z_{n}$ for $\phi$ tends to infinity
non-tangentially. Thus, given $\Om$ as in Lemma \ref{lem:pom} (a),
every orbit is eventually in $\Om$, so the grand orbits for $\phi$ are
in one-to-one correspondence with the forward orbits in $\Om$.  
Since $g$ is one-to-one on
$\Om$, we have that $\Si_{g}$ is one-to-one.
Also, since $g (\Om)$ has an inner-tangent at infinity, we
also have that $\Si_{g}$ is onto the space of grand orbits for $\alpha$.

In the parabolic non-zero-step case,
suppose that $\langle z\rangle\neq \langle w\rangle $. Then, referring
to Lemma \ref{lem:pom} (b),
$z_{n},w_{n}\in \Om (\rho ,m_{\rho})$ for some $\rho$ and for $n$ sufficiently
large. So $\Si_{g} (\langle z\rangle)\neq \Si_{g} (\langle w\rangle)$.
On the other hand if an orbit of $\alpha$ is given, $\{\zeta +nb \}$,
then $\zeta +nb\in g (U)$ for $n$
large. So $\Si_{g}$ is onto.
\end{proof}

\section{Main results}\label{sec:main}
Here we state our results regarding the classical
semiconjugations in the hyperbolic and parabolic non-zero-step
cases, we abbreviate these two cases by saying ``non-zero-step type'',
see Remark \ref{rem:step}. The
other parabolic case will be discussed later in the paper.
 
First we need a result similar to Lemma 2.61 of \cite{cowen-maccluer:2003tjm}.

\begin{proposition}\label{prop:main2}
Suppose $\phi$ is a self-map of $\Hh$ of non-zero-step type, as in
(\ref{eq:phi}) with $A\geq 1$. Let $(g,\alpha)$ be the semiconjugation
obtained in Theorem~A.
Given a solution $(\sigma,\tau)\in\cF
(\phi)$,
there is an analytic self-map $F$ of $\D$ such that
\begin{equation}\label{eq:factorh}
\sigma = F \circ g
\end{equation}
and
\begin{equation}\label{eq:automh}
F\circ \alpha =\tau \circ F.
\end{equation}
Moreover, $\sigma$ is canonical if and only if $F$ is.
\end{proposition}
Proposition \ref{prop:main2} says that in order to solve the Forward
Equation (\ref{eq:forwardeq}), it is enough to consider tha pair
$(g,\alpha)$ of Theorem A, and then solve (\ref{eq:forwardeq}) with
the self-map $\phi$ replaced by the automorphism $\alpha$. 
In particular, $F$ is a self-map of $\D$
which intertwines two automorphisms $\alpha$ and $\tau$.
This fact is helpful in determining the canonical solutions to
(\ref{eq:forwardeq}), because it is easier to
determine when $F$ is canonical (this will be done below in Section
\ref{sec:aut}). As a result we will obtain the following
characterizations of canonical solutions and thus answer Question
\ref{quest:uniq}. 

\begin{theorem}\label{thm:solh}
Assume $\phi$ is a self-map of
$\Hh$ of non-zero-step type. Let $(g,\alpha)$, $A$, $b$ be as in
Theorem~A, and let $\beta$ be
as in Remark \ref{rem:conj}.  
Given a solution $(\sigma,\tau)\in\cF
(\phi)$ where $\tau$ is not elliptic, let $\tilde{\beta}$ be the
automorphism of $\Hh$ that 
conjugates $\tau$ to its standard form. Then 
the following are equivalent:
\begin{enumerate}
\item [(a)] When $A>1$, 
$\tilde{\beta}^{-1}\circ \tau \circ \tilde{\beta} (z)= Az$, and there is 
$c>0$ so that $\si=\tilde{\beta}\circ (c\beta^{-1})\circ
g$. 

When $A=1$, 
$\tilde{\beta}^{-1}\circ \tau \circ \tilde{\beta} (z)= z+1$, and 
there is $d\in \R$ such that
$\si =\tilde{\beta} (\beta^{-1}\circ g +d)$. 
\item [(b)] $\sigma$ is canonical (see Definition \ref{def:canonical}).
\item [(c)] $\tilde{\beta}^{-1}\circ \si$ has the same univalence and
covering properties of 
$g$ in Lemma \ref{lem:pom}.
\item [(d)] There exists a pair of points $z,w\in\Hh$ for which
equality holds in (\ref{eq:maximal}),  not identically zero.
\end{enumerate}
\end{theorem}

\begin{corollary}\label{cor:uniqh}
Suppose $\phi$ is a self-map of $\Hh$ of non-zero-step  type,
written as in 
(\ref{eq:phi}) with $A\geq 1$. Consider the semiconjugations
$(g,\alpha (z)=Az+b)$ 
and $(\tilde{g},\tilde{\alpha} (z)=Az+\tilde{b})$,
obtained in Theorem A by renormalizing $\phi_{n}$ using
$\{\phi_{n} (z_{0})\}$ and $\{\phi_{n} (\tilde{z}_{0})\}$ respectively.
Then, when $A>1$,
\[
\frac{\tilde{g}+\frac{\tilde{b}}{A-1}}{i+\frac{\tilde{b}}{A-1}}
=\frac{g+\frac{b}{A-1}}{ g(\tilde{z}_0)+\frac{b}{A-1}};
\]
while, when $A=1$,
\[
\frac{\tilde{g}-i}{|\tilde{b}|}=\frac{g-g(\tilde{z}_0)}{|b|}.
\]
\end{corollary}

\section{The forward equation}\label{sec:forward}

We give the proof of Proposition \ref{prop:main2} here and we postpone the
proofs of Theorem \ref{thm:solh} and 
Corollary \ref{cor:uniqh} to a later section.

\begin{proof}[Proof of Proposition \ref{prop:main2}]
Assume first that $\phi$ is of hyperbolic type.
Let $\Om$ be the simply connected domain in $\Hh$ with inner-tangent at
infinity, such that $g$ is one-to-one on $\Om$, which is given by
Lemma \ref{lem:pom} (a). Recall that $g (\Om)$ also
has an inner-tangent at infinity, and that $\phi
\circ(g\mid \Om)^{-1} (w)= (g\mid \Om)^{-1} (\alpha (w))$ whenever
$w,\alpha (w)\in g (\Om)$.

Given $(\sigma ,\tau)\in \cF (\phi)$,
define $F (w)=\sigma \circ (g\mid \Om)^{-1}$ on $g (\Om)$. 
Given $w\in \Hh\setminus \Om$, let
$n_{0}=n_{0} (w)$ be the smallest integer such that the iterates
$\alpha_{n}(w)$ are in $\Om$ for $n\geq n_{0}$. Then, 
\begin{eqnarray*}
\tau_{n_{0}+k}^{-1}\circ F\circ \alpha_{n_{0}+k}
(w) & = & \tau_{n_{0}+k}^{-1}\circ \sigma \circ
(g\mid \Om)^{-1} \circ \alpha_{k}\circ \alpha_{n_{0}}(w)\\
& = &
\tau_{n_{0}}^{-1}\circ\tau_{k}^{-1}\circ  \sigma \circ \phi_{k}\circ
(g\mid \Om)^{-1} 
\circ  \alpha_{n_{0}}(w)\\
& = &
\tau_{n_{0}}^{-1}\circ \sigma \circ(g\mid \Om)^{-1}
\circ \alpha_{n_{0}}^{-1} (w)\\
& = &
\tau_{n_{0}}^{-1}\circ F \circ \alpha_{n_{0}}^{-1} (w)
\end{eqnarray*}
Therefore, $F$ is well-defined on all of $\Hh$ via the formula $F
(w)=\tau_{n_{0}}^{-1}\circ \sigma \circ(g\mid \Om)^{-1}
(\alpha_{n_{0}} (w))$, with the understanding that $n_{0}=0$ if
$w\in \Om$. Moreover, $F\circ \alpha =\tau \circ F$, so $F\in \cA
(\alpha ,\tau)$.

Now assume that $\phi$ is of parabolic non-zero-step type.
Let $(\sigma ,\tau)\in \cF (\phi)$, and let $w\in\Hh$.
Pick $\rho =2\rho_{\Hh} (w,i)$. By Lemma \ref{lem:pom} (b), 
for $n\geq n (w)$, $\alpha_{n}
(w)\in g (U)$, and 
$\left(g\mid U\right)^{-1} (\alpha_{n+k} (w))=\phi_{k} (
\left(g\mid U \right)^{-1}  (\alpha_{n} (w)))$, for 
$k=1,2,3,\dots$.
Thus, we define
\[
F (w)=\tau_{n}^{-1}\circ \sigma \circ \left(g\mid U
\right)^{-1}  \circ \alpha_{n} (w) 
\]
for some $n\geq n (w)$, and this definition does not depend on which
$n\geq n (w)$ is chosen. Then $F$ is an analytic self-map of $\Hh$, and 
$F\circ \alpha =\tau \circ F$. Moreover $F\circ g (z)=\sigma (z)$.

For the ``moreover part'', it is clear that $\sigma$ is canonical if
and only if both $F$ and $g$ are, so the claim follows from
Proposition \ref{prop:gcanon}. 
\end{proof}

The next section lists the possible intertwining maps between
automorphisms of the disk, and identify the ones that are canonical.

\section{Maps that intertwine automorphisms of the disk}\label{sec:aut}
In this section we study the forward equation (\ref{eq:forwardeq})
when $\phi :=\ga  \in \Aut (\D)$. Actually in this situation it does
not make sense to distinguish 
between forward and backward.

Since there are three cases for $\ga $ and three for $\tau$ this 
yields nine cases. However, by Remark \ref{rem:fpaut}, only seven may have
non-trivial solutions. We are also using (\ref{eq:conjphif}) and
(\ref{eq:conjtf}).
\subsection{The cases}\label{ssec:cases}
\begin{itemize}
\item   $\cA_{h\rightarrow h} (S,T) =\{\sigma:\Hh \rightarrow \Hh; T,S>1,
\mid 
\sigma (Sz) =T\sigma (z)\}$.
\item $\cA_{h\rightarrow p} (S,\pm) =\{\sigma:\Hh \rightarrow \Hh; S>1\mid
\sigma (Sz) =\sigma\pm 1\}$.
\item $\cA_{h\rightarrow e} (S,\theta) =\{\sigma:\Hh \rightarrow \D; S>1;
\theta\in [0,2\pi)\mid 
\sigma (Sz) =e^{i\theta}\sigma (z)\}$.
\item $\cA_{p\rightarrow p} (\pm,\pm) =\{\sigma:\Hh \rightarrow \Hh \mid
\sigma (z\pm 1) 
=\sigma (z)\pm 1\}$.
\item $\cA_{p\rightarrow h} (\pm,T) =\{\sigma:\Hh \rightarrow \Hh; T>1\mid
\sigma (z\pm 1) =T\sigma (z)\}$.
\item $\cA_{p\rightarrow e} (\pm,\theta) =\{\sigma:\Hh \rightarrow \D;
\theta\in [0,2\pi) 
\mid \sigma (z\pm 1)=e^{i\theta}\sigma\}$.
\item $\cA_{e\rightarrow e} (\varphi,\theta ) =\{\sigma:\D\rightarrow \D;
\theta,\varphi\in [0,2\pi)\mid
\sigma (e^{i\varphi}z) =e^{i\theta}\sigma (z)\}$.
\end{itemize}

For simplicity, we again restrict our attention to $\ga $ hyperbolic or
parabolic to start with. The case $\cA_{e\rightarrow e}$, when $\ga $
is elliptic, is actually
more complicated and will be postponed to a later section.

\begin{proposition}\label{pro:proposicio} We have
\begin{enumerate}
\item $\cA_{p\rightarrow p} (\pm,\pm) =\cA_{p\rightarrow p}
(+,+)\cup\cA_{p\rightarrow p} (-,-)$. Moreover, $\cA_{p\rightarrow p}
(+,+)=\cA_{p\rightarrow p} (-,-)$.
\item $\cA_{p\rightarrow h} (\pm,T)$ is empty.
\item If $1<S<T$, then  $\cA_{h\rightarrow h} (S,T)$ is empty; 
if $1<T<S$, then every $\sigma \in
\cA_{h\rightarrow h} (S,T)$ satisfies $\ntlim_{z\rightarrow
\infty}\sigma (z)/z=0$; while if $T=S$, then  $\cA_{h\rightarrow h}
(S,S)=\{cz \mid c>0\}$.
\end{enumerate}
\end{proposition}
In all three cases in this proposition, $\sigma$ is a self-map of
$\Hh$, so we can write 
\[
\sigma (z)=cz+p (z)
\]
with $0\leq c<\infty$, $\ima p (z)>0$ and $\ntlim_{z\rightarrow
\infty}p (z)/z=0$. 
We first observe that $\sigma  (z)/z$ has the {\sf
Lindel\"of property}, i.e., 
suppose $\ga  :[0,+\infty)\rightarrow \Hh$ is a path such that
$\lim_{t\rightarrow +\infty}\gamma (t)=\infty$, and suppose
\[
\lim_{t\rightarrow +\infty} \frac{\sigma  (\ga  (t))}{\ga  (t)}=a
\]
for some $a\in \C\cup\{\infty\}$. Then, 
\begin{equation}\label{eq:lindel}
a=\ntlim_{z\rightarrow
\infty}\frac{\sigma  (z)}{z}=c.
\end{equation}
In fact, since $p (z)/z$ misses the negative real axis,
composition with a conformal map to the disk will turn it into a
bounded function. So (\ref{eq:lindel})
follows from the usual Lindel\"of property of
bounded analytic functions in $\Hh$, see \cite{ahlfors1973} p.~40,
\begin{proof}[Proof of Proposition \ref{pro:proposicio}]
In (1), suppose that
$\sigma$ is a self-map of $\Hh$ 
such that $\sigma (z+1)=\sigma (z)-1$ (the other case is similar). 
Let $\ga_{n} (t)=i+nt$, $t\in [0,1]$, $n=0,1,2,3\dots$, and $\Ga
=\cup\ga_{n}$.  Then, since $\sigma (z+n)=\sigma (z)-n$, 
the function $\sigma (z)/z$ has limit $-1$ along
$\Ga$. Therefore, by the Lindel\"of's property, $\ntlim_{z\rightarrow
\infty}\sigma (z)/z=-1=c$, but $c\geq 0$.

The moreover part is trivial.

In (2), suppose that $\sigma (z+1)=T\sigma (z)$ with $T>1$. Let
$\Ga$ be the 
path defined above. Then, since $\sigma
(z+n)=T^{n}\sigma (z)$, the function $\sigma (z)/z$ tends to infinity
along $\Ga$. So,  by the Lindel\"of's property, $\ntlim_{z\rightarrow
\infty}\sigma (z)/z=\infty=c$, but $c<\infty$.

In (3), since $\sigma (S^{n}z)=T^{n}\sigma (z)$, we have 
\[
\frac{\sigma (S^{n}z)}{S^{n}z}=\left(\frac{T}{S} \right)^{n}\frac{\sigma
(z)}{z}=
c+\frac{p (S^{n}z)}{S^{n}z}\rightarrow c.
\] 
Thus $T\leq S$, and if $T<S$, then $c=0$. 
On the other hand, if $T=S$, then
$\sigma (z)/z= c$, which implies that $c>0$ and that $\sigma (z)=cz$.
\end{proof}

\begin{remark}\label{rem:fhe}
If $\sigma\in \cA_{h\rightarrow h} (S,T)$, then $\log
\sigma /\log T\in \cA_{h\rightarrow p} (S,+)$ (here $0<\ima \log z<\pi$).

Also, if $\sigma\in \cA_{h\rightarrow p} (S,\pm)$, then $e^{i\theta
\sigma }\in 
\cA_{h\rightarrow e} (S,\pm\theta)$. Conversely, if $\sigma \in
\cA_{h\rightarrow e} (S,\theta)$, and $\sigma \neq 0$, then $
(\log \sigma)/ (i\theta )\in \cA_{h\rightarrow p} (S,+)$, for any
choice of $\log$.   
Examples can be found in $\cA_{h\rightarrow e} (S,\theta)$ that do
have zeros.

Finally, if $\sigma\in \cA_{p\rightarrow p} (\pm,\pm)$, then 
$e^{i\theta\sigma }\in \cA_{p\rightarrow e} (\pm,\pm\theta)$.
Conversely, if $\sigma \in\cA_{p\rightarrow e} (+,\theta)$, and
$\sigma \neq 0$, then $(\log \sigma)/ (i\theta )\in \cA_{p\rightarrow
p} (+,+)$, for any choice of $\log$.   
Examples can be found in $\cA_{p\rightarrow e} (+,\theta)$ that do
have zeros. 
\end{remark}

\subsection{Canonical solutions}\label{ssec:canon}
We next describe all the possible canonical solutions in the set
$\cA$ which comprises all the cases in Section \ref{ssec:cases} above, except
$\cA_{e\rightarrow e}$.  
\begin{theorem}\label{thm:canon}
The only canonical solutions arising in $\cA$ are:  \begin{enumerate}
\item $\{cz \mid c>0\}=\cA_{h\rightarrow
h} (S,S)$;
\item $\{z+d \mid d\in \R\}\subset \cA_{p\rightarrow p}
(\pm,\pm)$.
\end{enumerate}
\end{theorem}
\begin{proof}
We first list the possible grand orbit spaces that arise depending on
the nature of $\ga$:
\begin{itemize}
\item $\cG_{S}=\Hh/\langle z\mapsto Sz \rangle$ ($S>1$);
\item $\cG_{+}=\Hh/\langle z\mapsto z+1 \rangle$;
\item $\cG_{\theta}= \D/\langle z\mapsto e^{i\theta}z \rangle$
($\theta \in (0,2\pi)$).
\end{itemize}
Notice that $\cG_{S}$ is a Riemann surface conformally equivalent to
an annulus, and $\cG_{+}$ is a Riemann surface conformally equivalent
to a punctured disk. Moreover $\cG_{S}$ is not conformally equivalent
to $\cG_{+}$, and if $1<T<S$, then $\cG_{S}$ is not conformally
equivalent to $\cG_{T}$.

A solution $\sigma$ in $\cA_{h\rightarrow h} (S,T)$ induces an analytic
map $\Sigma$ from $\cG_{S}$ to $\cG_{T}$. Imposing that $\sigma$ is
canonical, means that $\Sigma$ is a bijection, and since $\Sigma$ is
analytic, this automatically implies that $\Sigma$ is biholomorphic. Thus,
$\cG_{S}$ is conformally equivalent to $\cG_{T}$, and therefore $S=T$.

Likewise, a solution $\sigma$ in $\cA_{h\rightarrow p} (S,+)$ factors
down to an analytic map $\Sigma$ between 
$\cG_{S}$ and $\cG_{+}$, and since these two Riemann surfaces are not
conformally equivalent, none of these solutions $\sigma$ are canonical.

Finally, suppose $\sigma \in\cA_{p\rightarrow p} (+,+)$ induces an
automorphism on $\cG_{+}$. Then $\sigma \in \Aut (\Hh)$ and actually
$\sigma \in \Aut_{\infty} (\Hh)$. So $\sigma (z)=cz+d$, with $c>0$ and
$d\in \R$. However, since $\sigma (z+1)=\sigma (z)+1$, we must have
$c=1$. Hence, $\sigma (z)=z+d$
and these are all the possible canonical solutions 
that can be found in $\cA_{p\rightarrow p} (+,+)$.

The case of $\cG_{\theta}$ is not so straightforward.
Suppose $\sigma$ is in $\cA_{h\rightarrow e}
(S,\theta)$, or in $\cA_{p\rightarrow e} (\pm,\theta)$, and assume
that $\sigma$ is canonical. 
When $\theta/(2\pi)\in \Q$, we still obtain a Riemann surface
$\cG_{\theta}^{\star}= (\D\setminus \{0 \})/\langle z\mapsto
e^{i\theta}z \rangle$ which is conformally equivalent to a punctured
disk. Since $\sigma$ is canonical it must be
non-constant. Letting $Z=\sigma^{-1} (0)\subset\Hh$, we see that
$\sigma$ factors down to a conformal equivalence between the punctured
disk $\cG_{\theta}^{\star}$ and either $\cG_{S}\setminus P$ or
$\cG_{+}\setminus Q$, where $P$ and $Q$ are non-empty discrete sets of
punctures, corresponding to $Z$. Since no such equivalence exists,
$\sigma$ cannot be canonical.

Now assume that $\theta /(2\pi)\not\in \Q$.
Then $\sigma (z)=\sigma (w)$ implies that
$z$ and $w$ generate the same grand orbit. Therefore, there is $k\in
\Z$ such that $w=S^{k}z$, or $w=z+k$. Thus, $\sigma
(z)=e^{ik\theta}\sigma (z)$, so $\sigma (z)=\sigma
(w)=0$. Then $\sigma$ vanishes on
an infinite discrete set $Z$, yet is one-to-one on $\Hh \setminus Z$,
and this contradicts the argument principle.
\end{proof}

\section{Canonical semiconjugations for self-maps with non-zero step
}\label{sec:general}  

In this section we prove Theorem \ref{thm:solh} and Corollary
\ref{cor:uniqh}. 

\begin{proof}[Proof of Theorem \ref{thm:solh}]
The implication $(a)\Rightarrow (b)$ follows from Proposition
\ref{prop:gcanon},  since $\si$ is equal to $g$
post-composed by an automorphism of $\Hh$.

Conversely, by Proposition \ref{prop:main2} any solution $\sigma$ can be
written as $\sigma =F\circ g$, with $F\circ \alpha =\tau \circ F$, and
$\sigma$ is canonical if and only if $F$ is. 
By (\ref{eq:conjphif}) and (\ref{eq:conjtf}) 
\begin{equation}\label{eq:conj}
(F,\tau)\in \cF (\alpha)\Longleftrightarrow (\tilde{\beta}^{-1} \circ
F\circ \beta 
,\tilde{\beta}^{-1}\circ \tau \circ \tilde{\beta} )\in  \cF
(\beta^{-1}\circ \alpha 
\circ \beta)
\end{equation}
Choose $\tilde{\beta} $ so that $\tilde{\beta}^{-1}\circ \tau \circ
\tilde{\beta} $ is in 
standard form as in Section \ref{ssec:class}. 
Notice that $\tilde{\beta} ^{-1} \circ F\circ\beta  $ is canonical, and
when $A>1$ it intertwines $z\mapsto Az$ with $\tilde{\beta}^{-1}\circ \tau
\circ \tilde{\beta}$,
which is some other automorphism of $\Hh$ in
standard form. Therefore  
Theorem
\ref{thm:canon} implies that $\tilde{\beta} ^{-1} \circ F\circ\beta=cz$ for
some $c>0$, and also that
$\tilde{\beta}^{-1}\circ \tau \circ \tilde{\beta}  (z)=Az$. So, $F
(z)=\tilde{\beta} (c\beta^{-1} (z))$. In particular,
$\tilde{\beta}^{-1}\circ \sigma =c 
\beta^{-1}\circ g$.
On the other hand, when $A=1$, we assume without loss of generality
that $b>0$. Then $\tilde{\beta} ^{-1} \circ F\circ\beta$
intertwines $z\mapsto z+1$ with $\tilde{\beta}^{-1}\circ \tau
\circ \tilde{\beta} $. So by Theorem
\ref{thm:canon}, $\tilde{\beta} ^{-1} \circ F\circ\beta=z+d$ for some $d\in
\R$, and also that 
$\tilde{\beta}^{-1}\circ \tau \circ \tilde{\beta}  (z)=z+1$. 
So, $F(z)=\tilde{\beta} (\beta^{-1} (z)+d)$. In particular,
$\tilde{\beta}^{-1}\circ \sigma = 
\beta^{-1}\circ g+d$.
Hence, we have shown $(b)\Rightarrow (a)$.

Note that $(a)\Rightarrow (c)$ is clear from Lemma
\ref{lem:pom}. Conversely, by Proposition \ref{prop:main2} we can write
$\sigma =F\circ g$, 
with $F\circ \alpha =\tau \circ F$. We claim that
$F$ must be an automorphism of 
$\Hh$. This implies that $\sigma$ is canonical, so $(c)\Rightarrow (a)$.
To see  that $F$ is one-to-one, suppose $F (z)=F (w)$. Then 
\[
F (\alpha_{n} (z))=\tau_{n}\circ F (z)=\tau_{n}\circ F (w)=F (\alpha_{n} (w))
\]
Find $\rho >0$ large enough so that $z,w\in \De (i,\rho)$.
By the proof of the 'furthermore' part in Lemma \ref{lem:pom}, we can
choose $N$ large enough so that $\alpha_{n} (z),\alpha_{n} (w)\in g
(\De (z_{n},\rho))$, and so that $g$ is univalent on $\De
(z_{n},\rho)$, for all $n\geq N$. By hypothesis, we can find $n$ large
enough so that $\sigma =F\circ g$ is univalent on $\De
(z_{n},\rho)$.
Thus $\alpha_{n} (z)=\alpha_{n} (w)$, and $z=w$.

To see that $F$ is onto, recall that, by hypothesis,
$\tilde{\beta}^{-1}\circ\si$ has certain univalence and covering
properties at $\infty$. Thus,  $\tilde{\beta}^{-1}\circ F$ must  have
certain corresponding properties. Hence we can define
\[
G(z):=\alpha_n^{-1}\tilde{F}^{-1}(\tilde{\tau}_n^{-1}(z))
\] 
to be a global analytic inverse of $\tilde{F}$ in $\Hh$, where
$\tilde{\tau}=\tilde{\beta}^{-1}\circ\tau\circ\tilde{\beta}$, and $n$
is chosen large enough depending on $z$. 

Finally, $(a)\Rightarrow (d)$ is clear from (\ref{eq:maximal}). Conversely,
write $\sigma =F\circ g$. Then, 
\[
\rho_{\Hh} (F (g (z)),F (g (w)))=\rho_{\Hh} (g (z),g (w))\neq 0.
\]
So, by Schwarz's Lemma, $F$ is an automorphism of $\Hh$. Hence, $\sigma$
is canonical, and $(d)\Rightarrow (a)$.
\end{proof}

\begin{proof}[Proof of Corollary \ref{cor:uniqh}]
Let $\sigma =\tilde{g}$ and $\tau=\tilde{\alpha}$, then apply Theorem
\ref{thm:solh} and the fact that $\tilde{g}$ is canonical. Then plug
in $z=\tilde{z}_0$.
\end{proof}

\section{The parabolic zero-step case}\label{ssec:parabzs}
Proposition \ref{pro:proposicion} shows that when $\phi$ is parabolic 
zero-step, then $\cF (\phi)$ is empty. So instead of looking for pairs
$(\sigma ,\tau)$ where $\sigma$ is a self-map of $\Hh$ and $\tau$ is
an automorphism of $\Hh$, one is led to try to find pairs where
$\sigma$ is analytic and $\tau$ is linear, namely  solutions
to the Planar Equation.
\begin{definition}\label{def:pphi}
We write $\cP (\phi)$ for the family of all non-trivial solution pairs $(\sigma
,\tau)$ to (\ref{eq:planareq}), i.e., we ask that $\sigma$ be
non-constant. 
\end{definition}
Note that (\ref{eq:conjphif}) and (\ref{eq:conjtf}) still hold with
$\cF$ replaced by $\cP$ and $\alpha ,\beta \in \Aut (\C)$. So 
$\tau$ can always be put in standard form:
namely, either $\tau (z)=az$ for some $a\in \C\setminus \{0\}$,  or
$\tau (z)=z+1$. Again we will write
\begin{itemize}
\item $\cP_{e}(\phi,a)=\{\sigma:\Hh\rightarrow \C ;a\in \C\setminus \{0\}\mid
\sigma \circ \phi =a\sigma\}$.
\item $\cP_{p}(\phi,+)=\{\sigma:\Hh \rightarrow \C \mid \sigma \circ \phi
=\sigma+1\}$.
\end{itemize}
When $\phi$ is parabolic zero-step,
existence of 
solutions in $\cP_{p}(\phi,+)$ is found in 
\cite{baker-pommerenke:1979jlms}. In 
\cite{cowen:1981tams}
this dichotomy between disk
model and plane model is explained in the framework of the ``type
problem'' for simply-connected non-compact Riemann surfaces and the
uniformization theorem.

Fix a self-map of $\Hh$ of parabolic zero-step type, which can
therefore be written as
\begin{equation}\label{eq:phizs}
\phi (z) = z+p (z)
\end{equation}
$\ima p (z)>0$ and $\ntlim_{z\rightarrow \infty}p (z)/z=0$
(see the classification in Section \ref{ssec:class}).
Write $z_{n}=\phi_{n} (z_{0})=x_{n}+iy_{n}$, and let 
\[
\mu_{n} (z):=\frac{z-z_{n}}{z_{n+1}-z_{n}}\in \Aut (\C)
\]
be the automorphism of $\C$ which sends $z_{n}$ to $0$ and $z_{n+1}$
to $1$. Then consider the normalized iterates
\[
h_{n}=\mu_{n}\circ \phi_{n}.
\]
\begin{theoremB}[Theorem 1 of \cite{baker-pommerenke:1979jlms}]
The limit $h=\lim_{n\rightarrow \infty}h_{n}$ exists locally uniformly
in $\Hh$,
satisfies $h(z_{0})=0$, and solves
\[
h\circ \phi =h+1.
\]
\end{theoremB}
\begin{remark}\label{rem:parbzsorbit}
The zero-step condition implies that
\[
\frac{y_{n+1}}{y_{n}}\longrightarrow 1 \qquad \mbox{and}\qquad
\frac{x_{n+1}-x_{n}}{y_{n}}\longrightarrow 0,
\]
and again $y_{n+1}>y_{n}$ so $y_n\rightarrow L_{\infty}$.
Using the zero-step hypothesis one sees that $L_{\infty}= +\infty$.
\end{remark}
We now describe the univalence and covering properties of the
semiconjugation $h$ 
in the spirit of Section \ref{ssec:classuniv}.
\begin{lemma}\label{lem:bapom}
The sequence $\psi_{n}=h\circ \mu_{n}^{-1}-n$ converges uniformly on compact
subsets of $\C$ to the identity map on $\C$. 
Moreover, sequences $R_{k}, M_{k}\uparrow+\infty$ can be chosen so that
$h$ is univalent on the set 
\begin{equation}\label{eq:union}
U=\bigcup_{k=1}^{\infty}\bigcup_{n=M_{k}}^{\infty} (z_{n}+
(z_{n+1}-z_{n})R_{k}\D)
\end{equation}
Furthermore, the region $V=h
(U)$ has a {\sf
lateral-tangent at $+\infty$ with respect to $\C$}, i.e., given
$R_{n}\uparrow +\infty$ 
there is $t_{n}\uparrow+\infty$ such that $\{\rea z>t_{n},\
|\ima z|<R_{n}\}\subset V$.
\end{lemma}
\begin{proof}
First recall part (c)  of the proof of the main theorem in
\cite{baker-pommerenke:1979jlms}, where it shown that there is $m$ so
that $h (z_{n}+y_{n}s)$ is univalent in $|s|<1/5$ for all $n\geq m$.
Since $\mu_{n}^{-1} (\zeta)=z_{n}+\zeta (z_{n+1}-z_{n})$, and since by
Remark \ref{rem:parbzsorbit} 
\[
\frac{z_{n+1}-z_{n}}{y_{n}}=\frac{x_{n+1}-x_{n}}{y_{n}}+
i\left(\frac{y_{n+1}}{y_{n}}-1\right)\rightarrow 0 
\]
we see that given $R>1$ there is $N=N (R)\in \N$ such that
\[
\psi_{n} (\zeta) =h\circ \mu_{n}^{-1} (\zeta ) -n=h (z_{n}+\zeta
(z_{n+1}-z_{n}))-n
\]
is well defined and univalent for all $|\zeta |<R$ and all $n\geq N$. 

Moreover, $\psi_{n} (0)=h (z_{n})-n=0$ and $\psi_{n} (1)=h
(z_{n+1})-n=1$. So the sequence $\{\psi_{n} \}$ is normal on $\{|\zeta
|<R \}\setminus  \{0,1 \}$, and hence, by a standard argument, on
$\{|\zeta |<R \}$ as well. Now,
\[
\psi_n (h_{n})=h\circ \mu_{n}^{-1}\circ \mu\circ \phi_{n}-n= h
\]
So any normal sublimit $\psi$ of the $\psi_{n}$ is the identity map on
$h (\Hh)$. But since $h$ is an open map, this implies that $\psi_{n}$
 converges to the identity on $\C$.

To prove the univalence statement we proceed as in the proof of
Theorem 3 in \cite{pommerenke:1979jlms}. First note that, since
$\psi_{n}$ tends to the identity map uniformly on compact sets, for every
$R>1$ there is an integer $k_{R}$ such that
\begin{equation}\label{eq:overlap}
\left(\psi_{n} (R\D)\right)\cap \left(\psi_{n+k_{R}} (R\D)+
k_{R} \right) =\emptyset 
\end{equation}
for $n=1,2,3,\dots$

Now suppose that $h (s_{1})=h (s_{2})$ at points $s_{1}=\mu_n^{-1}
(\zeta_{1})$ and $s_{2}=\mu_{n+j}^{-1}
(\zeta_{2})$  with $j\geq 0$ and $|\zeta_{l}|<R$ for $l=1,2$. Then
\begin{eqnarray*}
\psi_n (\zeta_{1}) & = & h\circ \mu_n^{-1}
(\zeta_{1})-n= h (s_{1})-n\\
 & = & h (s_{2})-n=h\circ \mu_{n+j}^{-1}
(\zeta_{2})- (n+j)+j=\psi_{n+j} (\zeta_{2})+j
\end{eqnarray*}
So by (\ref{eq:overlap}) we must have $j\leq k_{R}$.

On the other hand,
\[
s_{2} = z_{n+j}+ (z_{n+j+1}-z_{n+j})\zeta_{2} =z_{n}+ (z_{n+1}-z_{n})t_{n}
\]
where 
\[
t_{n}=\frac{(z_{n+j}-z_{n})+
(z_{n+j+1}-z_{n+j})\zeta_{2}}{z_{n+1}-z_{n}}\rightarrow j-\zeta_{2} 
\]
as $n\rightarrow \infty$. So $s_{2},s_{1}\in \mu_{n}^{-1} (R'\D)$ with
$R'=R+k_{R}+1$. However, since $\psi_{n}$ converges to the identity
map on compact sets, there is $m_{R}$ such that $h$ is one-to-one on
$\mu_{n}^{-1} (R'\D)$ for all $n\geq m_{R}$. Therefore, we must have $n<m_{R}$.
In other words, $h$ is one-to-one on $\cup_{n\geq m_{R}}\mu_{n}^{-1}
(R\D)$.
Thus, given a
ball $R\D$, there is $m_{R }\in \N$ so that $h$ is univalent on the
set 
\[
U (R,m_{R}) =\bigcup_{n=m_R}^{\infty} (z_{n}+
(z_{n+1}-z_{n})R\D)=\bigcup_{n=m_R}^{\infty} \mu_{n}^{-1}(R\D) 
\]

The 'moreover' part is proved exactly as in the proof of Lemma
\ref{lem:pom}.

For the 'furthermore' part,
consider the half-strips $S_{t}=\{\rea z>t, |\ima z|<R/2 \}$.
Then by Lemma \ref{lem:bapom}, given $\epsilon <R/2$ there is $N$
such that, for all $n\geq 
N$, 
\[
|\psi_{n} (\zeta )-\zeta |=|h (z_{n}+\zeta (z_{n+1}-z_{n}))- (\zeta
+n) | <\epsilon  
\]
uniformly for $|\zeta |\leq R$. In particular, when $n>\max\{N, m_{R}
\}$, for every $|w|<R/2$,
$h\circ \mu_{n}^{-1} 
(\{|\zeta |=R \})$ winds around $w+n$ exactly once, i.e., $w+n\in h
(U)$. Therefore, $S_{t}\subset h (U)$ for $t$ large enough. 
\end{proof}

Next we prove the equivalent of Proposition \ref{prop:gcanon}.

\begin{proposition}\label{prop:hcanon}
The semiconjugation $h$ obtained in Theorem~B is canonical.
\end{proposition}
\begin{proof}
We need to show that the map $\Si_{h}$ which takes grand orbits of
$\phi$ to the orbits of $z\mapsto z+1$ on $\C$ is one-to-one and onto.
Again  write $\langle z\rangle$ for the grand orbit of $\phi$ generated by $z$.
Suppose that $\langle z\rangle\neq \langle w\rangle $. Write
$z_{n}=\phi_{n} (z)$ and $w_{n}=\phi_{n} (w)$. Then $h (z_{n})=h
(z)+n$ and $h (w_{n})=h (w)+n$. So, by Lemma \ref{lem:bapom},
there is $R$ large enough so that
for $n$ large enough $z_{n},w_{n}\in U(R,m_R)$. Since $z_{n+k}\neq w_{n}$
for all such $n$'s and $k\geq 0$, the univalence of $h$ implies that
$\Si_{h} (\langle z\rangle)\neq \Si_{h} (\langle w\rangle) $. On the
other hand given an orbit $\{a+n \}_{n=-\infty}^{+\infty}$ with $a\in
\C$, there is $N$ such that $a+n\in V=h (U)$ for $n\geq
N$. So $\Si_{h}$ is onto.
\end{proof}

We now state a result equivalent to Proposition \ref{prop:main2}.

\begin{proposition}\label{prop:parabnzs}
Suppose $\phi$ is a self-map of $\Hh$ of parabolic zero-step type. Let
$h$ be the semiconjugation obtained in Theorem~B. Given a solution
$(\sigma ,\tau)\in \cP (\phi)$, there is an entire function $F$ such
that
\[
\sigma =F\circ h
\] 
and 
\[
F (z+1)=\tau \circ F (z).
\]
Moreover, $\sigma$ is canonical if and only if $F$ is.
\end{proposition}

\begin{proof}
The proof is exactly the same as the one for Proposition \ref{prop:main2}.
\end{proof}

Next, we state a result equivalent to Theorem \ref{thm:solh}.

\begin{theorem}\label{thm:solp}
Assume that $\phi$ is a self-map of $\Hh$ of parabolic zero-step
type. Let $h$ be given as in Theorem~B. Given a solution $(\sigma
,\tau)\in \cP (\phi)$, let $\tilde{\beta}$ be the automorphism of $\C$
which conjugates $\tau$ to its standard form (either $z+1$ or $az$).
Then the following are equivalent.\begin{enumerate}
\item [(a)] There is $b\in \C$ such that
$\tau (z)=z+b$, $\tilde{\beta} (z)=bz$, and $\sigma =b(h+c)$ for some
$c\in \C$.
\item [(b)] $\sigma$ is canonical.
\item [(c)] $\tilde{\beta}^{-1}\circ\sigma$ has the same univalence
and covering properties as $h$ in Lemma \ref{lem:bapom}.
\end{enumerate}
\end{theorem}
During writing of this paper, we were informed that Contreras,
Diaz-Madrigal, and Pommerenke have also proved a statement similar to
(a)$\Leftrightarrow$ (c) above, but with a different univalence
requirement, see \cite{contreras}.
\begin{corollary}\label{cor:uniqp}
Suppose $\phi$ is a self-map of $\Hh$ of parabolic zero-step type.  Let
$h$ and $\tilde{h}$ be semiconjugations obtained as in Theorem~B, by
renormalizing $\phi_n$ using $\{\phi_n(z_0)\}$ and
$\{\phi_n(\tilde{z}_0)\}$ respectively. Then,
$\tilde{h}(z)=h(z)-h(\tilde{z}_0)$.
\end{corollary}

\subsection{Entire functions that intertwine automorphisms of the
plane}\label{ssec:entire}

This section parallels Section \ref{sec:aut}. 
Three cases arise.
\begin{itemize}
\item $\cE_{p\rightarrow p}=\{F \mbox{ entire }|F(z+1)=F(z)+1\}$;
\item $\cE_{p\rightarrow e}(a)=\{F \mbox{ entire }; a\in
\C\setminus\{0\}|F(z+1)=aF(z)\}$;
\item $\cE_{e\rightarrow e}(a,b)=\{F \mbox{ entire }; a,b\in
\C\setminus\{0\}|F(az)=bF(z)\}$.
\end{itemize}
Note that there cannot be $F$ entire such that $F(az)=F(z)+1$, because at $z=0$,
we would have $F(0)=F(0)+1$, which is impossible. So $\cE_{e\rightarrow p}$ is
empty. Also, we will not be interested in $\cE_{e\rightarrow e}(a,b)$ for the
moment, because our immediate concern is proving Theorem \ref{thm:solp}.       

Next we analyze the canonical maps.
\begin{theorem}\label{thm:canonp}
The only canonical maps in $\cE_{p\rightarrow p}$ are $F(z)=z+c$, $c\in \C$.
Moreover, $\cE_{p\rightarrow e}(a)$ does not contain any canonical maps.
\end{theorem}
\begin{proof}
Suppose $F\in \cE_{p\rightarrow p}$, then $F$ induces a map $\Sigma$ from
$\C\setminus\langle z\mapsto z+1\rangle$ to itself, which is analytic. Imposing
that $\Sigma$ be one-to-one and onto, means that $\Sigma$ is an automorphism of
$\C\setminus\langle z\mapsto z+1\rangle$ 
and therefore $F$ must also be an automorphism of $\C$. Thus, $F(z)=dz+c$ and
since it must also intertwine $z\mapsto z+1$ with itself, we have $d=1$.

For the 'Moreover' part, assume first that $|a|\neq 1$. Then, is a compact
Riemann surface. Hence there are no 
biholomorphisms between $\C\setminus\langle z\mapsto z+1\rangle$ and $\C\setminus
\langle z\mapsto az\rangle$. 
If $|a|=1$, assume first that $a=e^{i\theta}$ with $\theta/(2\pi)\in \Q$. Let
$Z=F^{-1}(0)\in \C$. Then $F$ canonical would induce a biholomorphism between a
punctured infinite cylinder and the punctured plane. Since there are no such
maps, the conclusion holds. Now assume that $\theta/(2\pi)\not\in
\Q$. Then $F(z)=F(w)$ and $F$ canonical imply that $z$
and $w$ generate the same grand orbit. So there is $n\in \N$, such that,
say, $w=z+n$. Hence $F(z)=e^{in\theta}F(z)$, so $F(z)=F(w)=0$. Then
$F$ vanishes on an infinite discrete set $Z$, yet is one-to-one on
$\C\setminus Z$, and this contradicts the argument principle.
\end{proof}

\subsection{Canonical semiconjugations for self-maps of parabolic
zero-step type}

\begin{proof}[Proof of Theorem \ref{thm:solp}]
The implication (a) $\Rightarrow$ (b) follows from Proposition
\ref{prop:hcanon}. Conversely, by Proposition \ref{prop:parabnzs}, 
any solution $\sigma$ can be
written as $\sigma =F\circ h$, with $F(z+1)=\tau \circ F(z)$, and
$\sigma$ is canonical if and only if $F$ is, if and only if
$\tilde{\beta}^{-1}\circ F$ is.
Moreover,
$\tilde{\beta}^{-1}\circ F$ intertwines $z\mapsto z+1$ with either 
$z\mapsto z+1$ or $z\mapsto az$ ($a\neq 0$). Hence, by Theorem
\ref{thm:canonp}, $\tilde{\beta}^{-1}\circ F$ is canonical implies 
$\tilde{\beta}^{-1}\circ F(z)=z+c$, for some $c\in \C$, and
$\tilde{\beta}^{-1}\circ \tau\circ\tilde{\beta}(z)=z+1$.
So $\tau(z)=z+b$, $\tilde{\beta}(z)=bz$, and $F(z)=bz+bc$.
Therefore, (b) $\Rightarrow$ (a). 

The implication (a) $\Rightarrow$ (c) follows from Lemma
\ref{lem:bapom}.
Conversely,  by Proposition \ref{prop:parabnzs}, we can write $\sigma
=F\circ h$, with $F(z+1)=\tau \circ F(z)$.
Using the covering and univalence properties we then show that $F$
must be an automorphism of $\C$ in exactly the same way as in the
proof of Theorem \ref{thm:solh} (c) $\Rightarrow$ (a).
\end{proof}

\begin{proof}[Proof of Corollary \ref{cor:uniqp}]
Let $\sigma =\tilde{g}$, then apply Theorem
\ref{thm:solp} and the fact that $\tilde{g}$ is canonical. Then plug
in $z=\tilde{z}_0$.
\end{proof}
%

\begin{thebibliography}{BPC03}

\bibitem[Ahl73]{ahlfors1973}
Lars~V. Ahlfors.
\newblock {\em Conformal invariants: topics in geometric function theory}.
\newblock McGraw-Hill Book Co., New York, 1973.
\newblock McGraw-Hill Series in Higher Mathematics.

\bibitem[BP79]{baker-pommerenke:1979jlms}
I.~N. Baker and Ch. Pommerenke.
\newblock On the iteration of analytic functions in a halfplane. {II}.
\newblock {\em J. London Math. Soc. (2)}, 20(2):255--258, 1979.

\bibitem[BPC03]{bracci-pc:2003jyv}
Filippo Bracci and Pietro Poggi-Corradini.
\newblock On {V}aliron's theorem.
\newblock In {\em Future trends in geometric function theory, RNC Workshop
  Jyv\"{a}skyl\"{a}}, volume~92, pages 39--55. University of Jyv\"{a}skyl\"{a},
  Department of Mathematics and Statistics, 2003.

\bibitem[BS97]{bourdon-shapiro:1997mams}
Paul~S. Bourdon and Joel~H. Shapiro.
\newblock Cyclic phenomena for composition operators.
\newblock {\em Mem. Amer. Math. Soc.}, 125(596):x+105, 1997.

\bibitem[CDMP]{contreras}
M.~Contreras, S.~Diaz-Madrigal, and Ch. Pommerenke.
\newblock {Some remarks on Abel equation in the unit disk}.
\newblock Preprint.

\bibitem[CM03]{cowen-maccluer:2003tjm}
Carl~C. Cowen and Barbara~D. MacCluer.
\newblock Schroeder's equation in several variables.
\newblock {\em Taiwanese J. Math.}, 7(1):129--154, 2003.

\bibitem[Cow81]{cowen:1981tams}
Carl~C. Cowen.
\newblock Iteration and the solution of functional equations for functions
  analytic in the unit disk.
\newblock {\em Trans. Amer. Math. Soc.}, 265(1):69--95, 1981.

\bibitem[PC00]{pc:2000finn}
Pietro Poggi-Corradini.
\newblock Canonical conjugations at fixed points other than the
  {D}enjoy-{W}olff point.
\newblock {\em Ann. Acad. Sci. Fenn. Math.}, 25(2):487--499, 2000.

\bibitem[PC03]{pc:2003iberoam}
Pietro Poggi-Corradini.
\newblock Backward-iteration sequences with bounded hyperbolic steps for
  analytic self-maps of the disk.
\newblock {\em Rev. Mat. Iberoamericana}, 19(3):943--970, 2003.

\bibitem[Pom79]{pommerenke:1979jlms}
Ch. Pommerenke.
\newblock On the iteration of analytic functions in a halfplane.
\newblock {\em J. London Math. Soc. (2)}, 19(3):439--447, 1979.

\bibitem[Val31]{valiron:1931bsm}
Georges Valiron.
\newblock Sur l'it\'eration des fonctions holomorphes dans un demi-plan.
\newblock {\em Bulletin Sc. math.}, 55(2):105--128, 1931.

\bibitem[Val54]{valiron1954}
Georges Valiron.
\newblock {\em Fonctions analytiques}.
\newblock Presses Universitaires de France, Paris, 1954.

\end{thebibliography}

\def\cprime{$'$}

\end{document}